\documentclass{article}

\usepackage[latin1]{inputenc}
\usepackage{subfigure}
\usepackage{multicol}

\usepackage{graphicx}

\usepackage{amsmath, xypic}

\newtheorem{thm}{Theorem}

\newtheorem{defn}{Definition}
\newtheorem{lem}{Lemma}
\newtheorem{prop}{Proposition}

\newtheorem{exm}{Example}

\newtheorem{rem}{Remark}

\newcommand{\la}{\lambda}

\newlength{\cellsz}
\newcounter{cellsize}
\newcommand{\setcellsize}[1]{%
  \setcounter{cellsize}{#1}%
  \setlength{\cellsz}{\value{cellsize}\unitlength}}%
\setcellsize{12}%
\newcommand\cellify[1]{\def\thearg{#1}\def\nothing{}%
\hbox to 0pt{{\begin{picture}(\value{cellsize},\value{cellsize})
  \put(0,0){\line(1,0){\value{cellsize}}}
  \put(0,0){\line(0,1){\value{cellsize}}}
  \put(\value{cellsize},0){\line(0,1){\value{cellsize}}}
  \put(0,\value{cellsize}){\line(1,0){\value{cellsize}}} \end{picture}
\hss}}
\vbox to \cellsz{ \vss \hbox to \cellsz{\hss$#1$\hss} \vss}}
\newcommand\tableau[1]{\vcenter{\vbox{\let\\\cr
\baselineskip -16000pt \lineskiplimit 16000pt \lineskip 0pt
\ialign{&\cellify{##}\cr#1\crcr}}}}
\newcommand\tabl[1]{\vtop{\let\\\cr
\baselineskip -16000pt \lineskiplimit 16000pt \lineskip 0pt
\ialign{&\cellify{##}\cr#1\crcr}}}
\newlength{\varcellsz}
\newcounter{varcellsize}
\newcommand{\setvarcellsize}[1]{%
  \setcounter{varcellsize}{#1}%
  \setlength{\varcellsz}{\value{varcellsize}\unitlength}}%
\setvarcellsize{10}%
\newcommand\varcellify[1]{\def\varthearg{#1}\def\varnothing{}%
\hbox to 0pt{{\begin{picture}(\value{varcellsize},\value{varcellsize})
  \put(0,0){\line(1,0){\value{varcellsize}}}
  \put(0,0){\line(0,1){\value{varcellsize}}}
  \put(\value{varcellsize},0){\line(0,1){\value{varcellsize}}}
  \put(0,\value{varcellsize}){\line(1,0){\value{varcellsize}}} \end{picture}
\hss}}
\vbox to \varcellsz{ \vss \hbox to \varcellsz{\hss$#1$\hss} \vss}}
\newcommand\vartableau[1]{\vcenter{\vbox{\let\\\cr
\baselineskip -16000pt \lineskiplimit 16000pt \lineskip 0pt
\ialign{&\varcellify{##}\cr#1\crcr}}}}
\newcommand\vartabl[1]{\vtop{\let\\\cr
\baselineskip -16000pt \lineskiplimit 16000pt \lineskip 0pt
\ialign{&\varcellify{##}\cr#1\crcr}}}

\author{Ekaterina A. Vassilieva}
\title{Long Cycle Factorizations : Bijective Computation in the General Case}

\begin{document}
\maketitle
\paragraph{Abstract.}
This paper is devoted to the computation of the number of ordered factorizations of a long cycle in the symmetric group where the number of factors is arbitrary and the cycle structure of the factors is given. Jackson (1988) derived the first closed form expression for the generating series of these numbers using the theory of the irreducible characters of the symmetric group. Thanks to a direct bijection we compute a similar formula and provide the first purely combinatorial evaluation of these generating series.

%


\section{Introduction}
\label{sec : intro}
For integer $n$ we note $S_n$ the symmetric group on $n$ elements and $\gamma_n$ the permutation in $S_n$ defined by $\gamma_n = (12\ldots n)$. If $r$ is an integer we call {\it strictly increasing subsequence of $1\ldots r$} any sequence of the form $(i_1, i_2, \ldots i_u)$ where $1\leq i_1 < i_2 < \ldots < i_u \leq r$. Given such a subsequence $t$ containing $i$, we define $succ_t(i)$ the index following $i$ in $t$. If no such index exists $succ_t(i)$ is the first index of the sequence or $i$ itself if $t=(i)$.\\
This paper is devoted to the computation of the numbers $k^n_{p_1, p_2,\ldots, p_r}$ of factorizations of $\gamma_n$ as an ordered product of permutations $\alpha_1\alpha_2\ldots \alpha_r = \gamma_n$ such that for $1\leq i\leq r$, $\alpha_i$ belongs to $S_n$ and is composed of exactly $p_i$ disjoint cycles. More precisely, we use a direct bijection to show the following formula:

\begin{thm}[Main result]
\label{thm : main}
\begin{align}
\label{eqmain} \frac{1} {(n-1)!^{r-1}}&\sum_{p_1, p_2,\ldots, p_r} {k^n_{p_1, p_2,\ldots, p_r}} \prod_{1\leq i \leq r}x_i^{p_i} =  \sum_{p_1, p_2,\ldots, p_r}\sum_{{\bf a}}\Delta_r({\bf a})\binom{n}{\bf a}\prod_{1\leq i \leq r}\binom{x_i}{p_i}
\end{align}
\noindent The last sum runs over sequences ${\bf a} = (a_t)$ of $2^r-1$ non-negative integers $a_t$ with the index $t$ being any non empty strictly increasing subsequence of integers of $1\ldots r$ such that
\begin{align}
\label{condmain} \sum_t{a_t} = n,\;\;\;\;\sum_{t;l\notin t}{a_t} = p_l \;\mbox{   for   }\; 2\leq l \leq r,\;\;\;\;\sum_{t;1\notin t}{a_t} = p_1-1
\end{align}
Furthermore, the multinomial coefficient is defined by $$\binom{n}{\bf a} = \frac{n!}{\prod_t{a_t!}}$$
Finally, $\Delta_r({\bf a})$ is the determinant of the $r\times r$ matrix with coefficients $m_{i,j}$, $1\leq i,j\leq r$, where
\begin{center}
$m_{i,i} = p_i$ ($1\leq i \leq r$),\\
$m_{i,i+1} = -p_{i+1} - \sum_{j\neq i} m_{j,i+1}$ ($1\leq i \leq r-1$),\\
$m_{r,1} = 1-p_{1} - \sum_{j\neq r} m_{j,1}$,\\
for $j \neq i+1$ (modulo $r$), $m_{i,j} = -\sum_{t;i\in t, succ_t(j)-j \geq j-i+1}a_t$\\
\end{center}
 (the subtractions  on indices are modulo $r$).
\end{thm}

Let $\lambda = (\lambda_1, \lambda_{2},...,\lambda_p) \vdash n$ an integer partition of $n$ with $\ell(\lambda) = p$ parts sorted in decreasing order. We note $C_\lambda$ the conjugacy class of $S_n$ containing the permutations of cycle type $\lambda$ and $m_\lambda({\bf x})$ and $p_\lambda({\bf x})$ the monomial and power sum symmetric functions respectively indexed by $\lambda$ on indeterminate ${\bf x}$. Given $r$ integer partitions $\lambda^1, \lambda^2,\ldots, \lambda^r$ of $n$, a more refined problem is to compute the numbers $k^n_{\lambda^1, \lambda^2,\ldots, \lambda^r}$ of ordered factorizations $\alpha_1\alpha_2\ldots \alpha_r$ of $\gamma_n$ such that for for $1\leq i\leq r$, $\alpha_i$ belongs to $C_{\lambda^i}$. As a corollary of Theorem \ref{thm : main} we have:

\begin{thm}[Corrolary]
\label{thm : corro}

\begin{align}
\nonumber \frac{1}{(n-1)!^{r-1}}&\sum_{\lambda^1, \lambda^2,\ldots, \lambda^r \vdash n} {k^n_{\lambda^1, \lambda^2,\ldots, \lambda^r}} \prod_{1\leq i \leq r}m_{\lambda^i}({\bf x}^i)\\
\label{eqcorr}& =  \sum_{\lambda^1, \lambda^2,\ldots, \lambda^r \vdash n}\frac{\sum_{{\bf a}}\Delta_r({\bf a})\binom{n}{\bf a}}{\prod_i \binom{n-1}{\ell(\lambda_i)-1}}\prod_{1\leq i \leq r}p_{\lambda^i}({\bf x}^i)
\end{align}

\noindent where $\ell(\lambda^i)$ is substituted to $p_i$ in the definition of ${\bf a}$ in (\ref{condmain}).
\end{thm}

\subsection{Background}

Despite the attention the problem received over the past twenty years no closed formulas are known for the coefficients $k^n_{p_{1},\ldots,p_{r}}$ and $k^n_{\la^{1},\ldots,\la^{r}}$ except for very special cases. Using characters of the symmetric group and a combinatorial development, Goupil and Schaeffer \cite{GS} derived an expression for $k^n_{\lambda^1,\lambda^2}$ ($r=2$) as a sum of positive terms. This work has been later generalized by Poulalhon and Schaeffer \cite{PS} and Irving \cite{JI} but, as a rule, the formulas obtained are rather complicated. Using the theory of the irreducible characters of the symmetric group, Jackson \cite{DMJ} computed an elegant expression for the generating series in the LHS of (\ref{eqmain}) for arbitrary $r$ and an arbitrary permutation of $S_n$ instead of $\gamma_n$. This later result shows that the coefficients in the expansion of this generating series in the basis of $\binom{x_i}{p_i}$ can be derived as closed form formulas but fails to provide a combinatorial interpretation. Schaeffer and Vassilieva in \cite{SV}, Vassilieva in \cite{V} and Morales and Vassilieva in \cite{MV} provided the first purely bijective computations of the generating series in (\ref{eqmain}) and (\ref{eqcorr}) for $r=2,3$. In a recent paper, Bernardi and Morales \cite{BM} addressed the problem of finding a general combinatorial proof of Jackson's formula for the factorizations of $\gamma_n$. Using an argument based on several successive bijections and a probabilistic puzzle, they provide a complete proof for the cases $r=2,3$ and a sketch for $r=4$.
In the present paper we generalize and put together all the ideas developed in our previous articles (\cite{SV}, \cite{V} and \cite{MV}) and make them work in the general context of $r$-factorizations of $\gamma_n$. We prove theorems (\ref{thm : main}) and (\ref{thm : corro}) thanks to a direct (single step) bijection. The combinatorial ingredients we use and the bijection itself are described in sections \ref{secing} and \ref{secbij}. Section (\ref{secproof}) proves that the bijection is indeed one-to-one. While (\ref{eqmain}) is similar to Jackson's formula in \cite{DMJ}, the two expressions are different. We address their equivalence in section \ref{seceq}.

\section{Cacti, partitioned cacti and cactus trees}
\label{secing}
\subsection{Cacti}
\label{subseccacti}
Factorizations of $\gamma_n$ can be represented as {\it $r$-cacti} (short {\em cacti}), i.e $2-$cell decompositions of an oriented surface of arbitrary genus into a finite number of vertices ($0-$cells), edges ($1-$cells) and faces ($2-$cells) homeomorphic to open discs, with $n$ black and one white face such that all the black faces are $r$-gons and not adjacent to each other. They are defined up to an homeomorphism of the surface that preserves its orientation, the type of cells and incidences in the graph.  We consider {\it rooted} cacti, i.e. cacti with a marked black face. We assume as well that within each $r$-gon the $r$ vertices are colored with $r$ distinct colors so that moving around the $r$-gons counter-clockwise the vertex of color $i+1$ (modulo $r$) follows the vertex of color $i$. Each black $r$-gon is labeled with an index in $\{1,2,\ldots,n\}$ such that the marked $r$-gon is labeled $1$ and that moving around the white face starting from the edge linking the vertex of color $1$ and the vertex of color $r$ in this marked $r$-gon, the $i$-th edge connecting a vertex of color $1$ and a vertex of color $r$ belongs to the black $r$-gon of index $i$. 
\begin{prop}[\cite{V}]
\label{propcac}
Cacti as defined above are in bijection with $r$-tuples of permutations $(\alpha_1, \alpha_2,\ldots, \alpha_r)$ such that $\alpha_1\alpha_2\ldots\alpha_r =\gamma_n$. Each vertex of color $i$ corresponds to a cycle of $\alpha_i$ defined by the sequence of the indices of the $r$-gons incident to this vertex. 
\end{prop}
 As a consequence, $r$-cacti with $p_i$ vertices of color $i$ are counted by $k^n_{p_1,p_2,\ldots,p_r}$. \\

\begin{figure}[htbp]
  \begin{center}
  \includegraphics[width=0.33\textwidth]{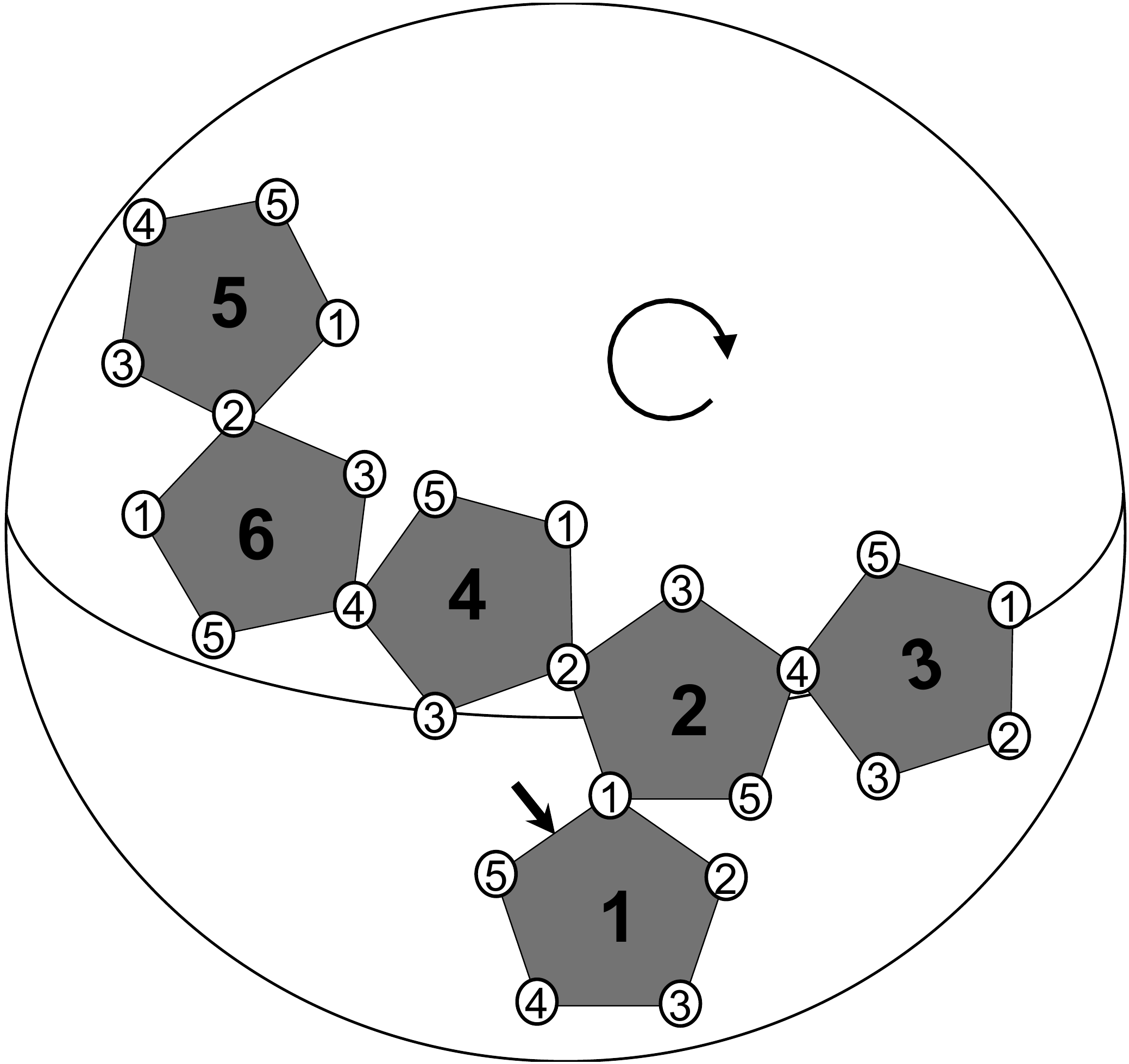}\hspace{10mm}
    \includegraphics[width=0.33\textwidth]{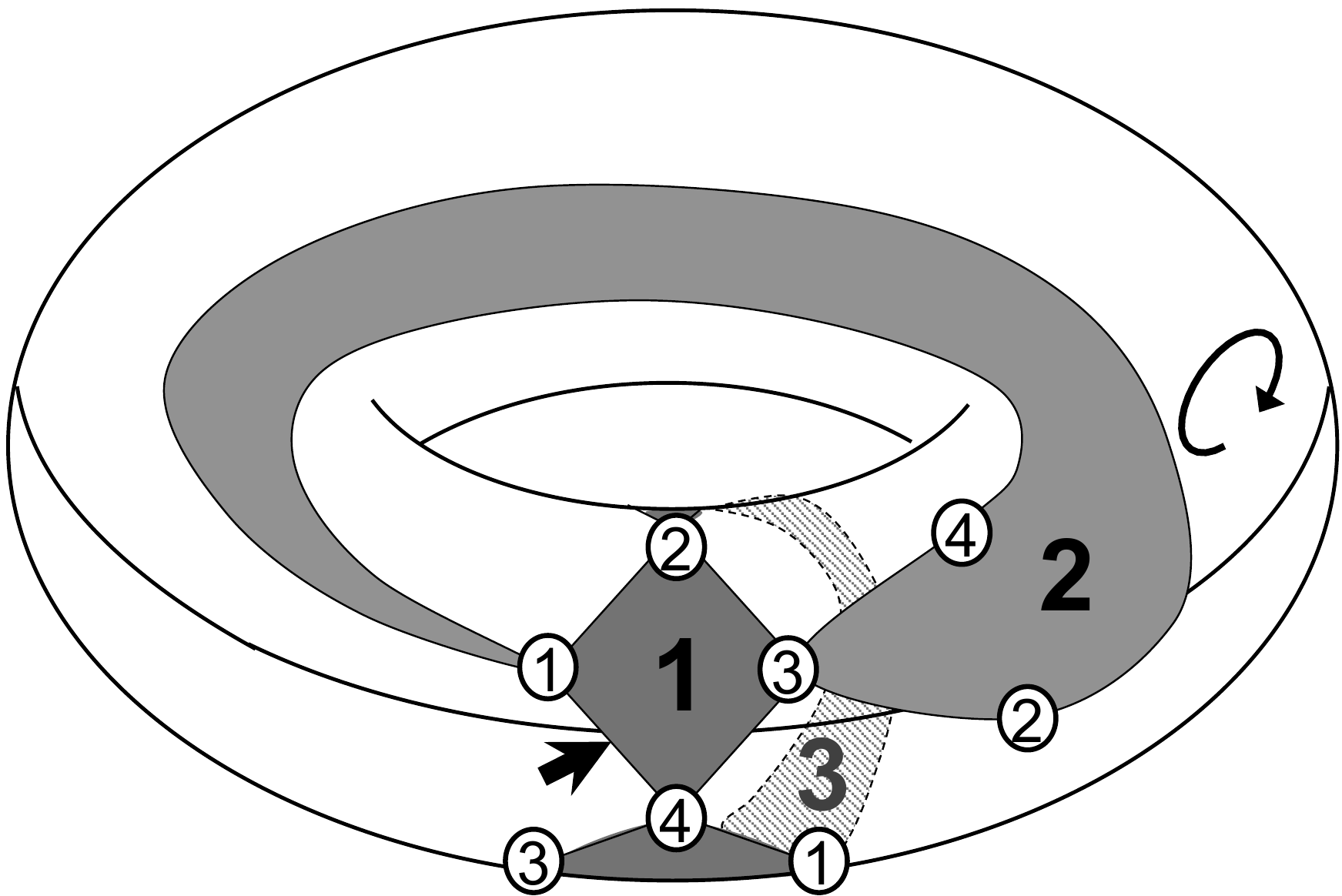}
    \caption{A $5$-cactus embedded on a surface of genus $0$ (left) and a $4$-cactus embedded on a surface of genus $1$ (right)}
     \label{cacti} \label{fig:example}
  \end{center}
\end{figure}

\begin{exm}
\label{exmcac}
Figure \ref{cacti} depicts a $5$-cactus corresponding to the factorizations of $\alpha_1\alpha_2\alpha_3\alpha_4\alpha_5 = \gamma_6$ with $\alpha_1 = (12)(3)(4)(5)(6)$, $\alpha_2 = (1)(24)(3)(56)$, $\alpha_3 = (1)(2)(3)(4)(5)(6)$, $\alpha_4 = (1)(23)(46)(5)$, $\alpha_5 = (1)(2)(3)(4)(5)(6)$ and a $4$-cactus described by $\alpha_1 = (12)(3)$, $\alpha_2 = (13)(2)$, $\alpha_3 = (12)(3)$, $\alpha_4 = (13)(2)$.
\end{exm}

\begin{rem}
\label{remcac}
Alternatively, a cactus can be seen as a set of $n$ $r$-tuples of integers such that each integer of $1\ldots n$ is used exactly once in the $i$-th positions of the $r$-tuples ($1\leq i \leq r$). Moving around the white face of the cactus according to the surface orientation and starting with the edge linking the vertices of color $1$ and $r$ of the root $r$-gon, we define a labeling of the edges. We assign label $i$ to the $i$-th edge linking vertices of color $1$ and $r$ during the traversal. The $j$-th edge linking vertices of color $r$ and $r-1$ during the same traversal is indexed by $j$ and so on for all the colors. The $n$ $r$-tuples defined by the edge labeling of the $n$ $r$-gons is an equivalent description of the initial cactus.     
\end{rem}

\begin{lem}
\label{lemcac}
One can show that if the label of an $r$-gon of a cacti described by $(\alpha_1, \alpha_2,\ldots, \alpha_r)$ is $i$ then the index (as defined in Remark \ref{remcac}) of the edge linking the vertices of color $l$ and $l-1$ is $\alpha_r^{-1}\alpha_{r-1}^{-1}\ldots\alpha_{l}^{-1}(i)$ for $2\leq l \leq r$ and $i$ for $l=1$.
\end{lem}

\begin{exm} The edge labeling defined above for the $5$-cactus of Example \ref{exmcac} is shown on Figure \ref{cactipart}.
\end{exm}

\subsection{Partitioned cacti}

Cacti are non planar non recursive objects intractable to compute directly in the general case. Our bijective construction relies on the use of {\it partitioned cacti} that we define as follows:

\begin{defn}[Partitioned cacti] Let $C^n(p_1,p_2,\ldots,p_r)$ be the set of couples composed of a cactus (as defined in Section \ref{subseccacti}) and a $r$-tuple of partitions $(\tilde{\pi}_1,\ldots,\tilde{\pi}_r)$ such that $\tilde{\pi}_i$ is a partition with $p_i$ blocks on the set of vertices of color $i$.
\end{defn}

\begin{rem}
\label{rempartcac}
Using Proposition \ref{propcac}, we state that partitioned cacti in $C^n(p_1,p_2,\ldots,p_r)$ are in bijection with the $2r$-tuples $(\alpha_1,\ldots,\alpha_r,\pi_1,\ldots,\pi_r)$ where $\alpha_1\ldots\alpha_r =\gamma_n$ and $\pi_i$ is a partition on the set of integers $1..n$ composed of exactly $p_i$ blocks stable by $\alpha_i$. (As such the blocks of $\pi_i$ are unions of cycles of $\alpha_i$). Two vertices $u$ and $v$ of color $i$ in the partitioned cactus belong to the same block if and only if the corresponding cycles of $\alpha_i$ belong to the same block of $\pi_i$.
\end{rem}

\begin{figure}[htbp]
  \begin{center}
 \includegraphics[width=0.45\textwidth]{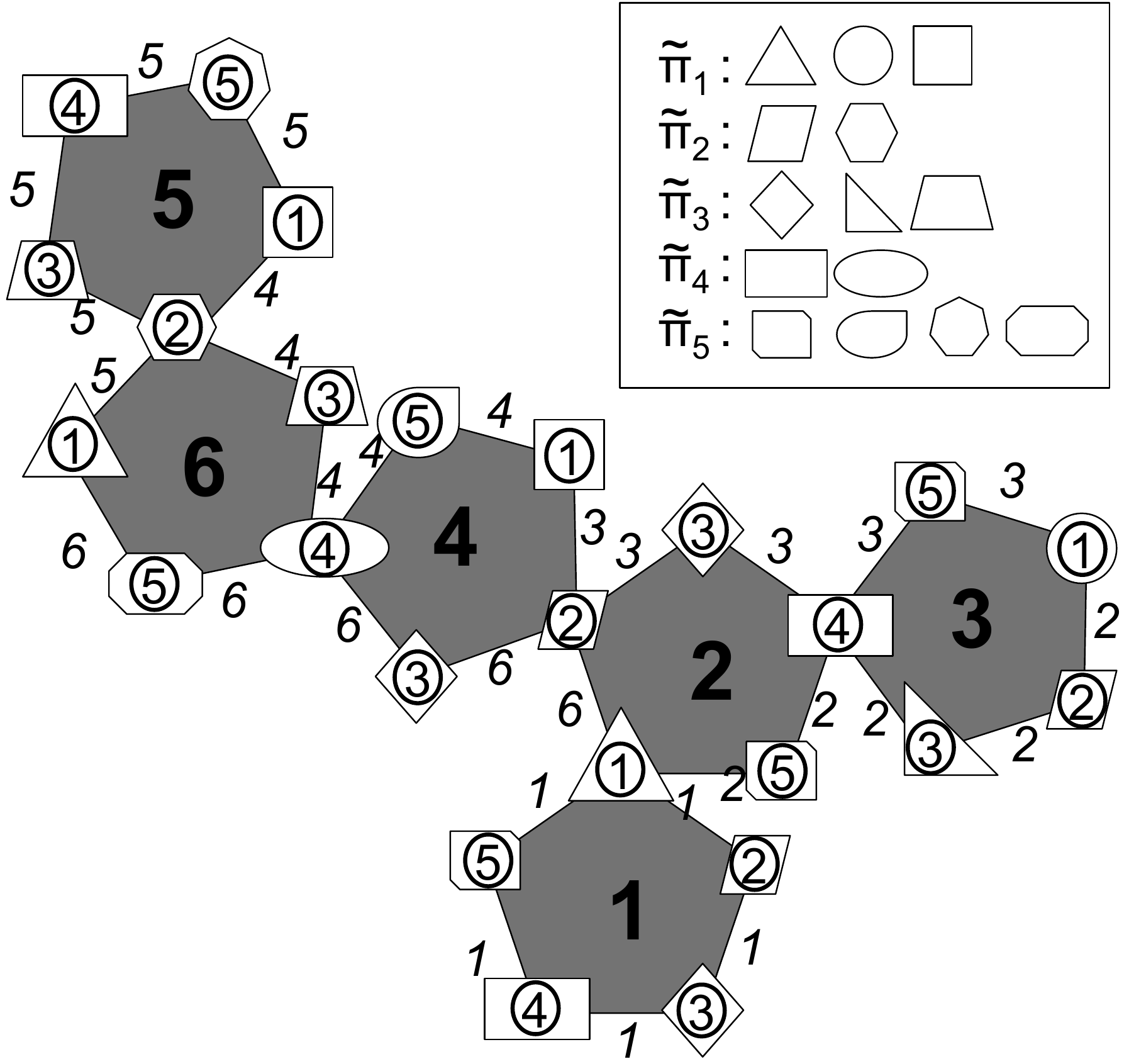}
    \caption{A partitioned cactus with the additional edge labeling defined in Remark \ref{remcac}.}
       \label{cactipart}
  \end{center}
\end{figure}

\begin{exm} We use geometric shapes in Figure \ref{cactipart} to represent an example of partitions on the sets of vertices of the $5$-cactus in Example \ref{exmcac}. The equivalent numeric set partitions are $\pi_1 = \{1,2,6\}\{3\}\{4,5\}$, $\pi_2 = \{1,2,3,4\}\{5,6\}$, $\pi_3 = \{1,2,4\}\{3\}\{5,6\}$, $\pi_4 = \{1,2,3,5\}\{4,6\}$, and $\pi_5 = \{1,2,3\}\{4\}\{5\}\{6\}$.
\end{exm}

Similarly to \cite{V}, the numbers $k^n_{p_1,p_2,\ldots,p_r}$ and the cardinalities $\mid C^n(p_1,p_2,\ldots,p_r)\mid$ are linked by the relation:

\begin{align}
\label{eqlink} \sum_{p_1, p_2,\ldots, p_r} {k^n_{p_1, p_2,\ldots, p_r}} \prod_{1\leq i \leq r}x_i^{p_i} =  \sum_{p_1, p_2,\ldots, p_r}\mid C^n(p_1,p_2,\ldots,p_r)\mid\prod_{1\leq i \leq r}({x_i})_{p_i}
\end{align}

\noindent where $({x})_{p} =x(x-1)\ldots(x-p+1)$. Partitioned cacti are actually one-to-one with decorated recursive tree structures that we define in the next section.

\subsection{Cactus trees}

We look at non classical tree-like structures with colored vertices and various types of children. More specifically we work with recursive non cyclic graphs rooted in a given vertex such that all the vertices are colored with $1,2,\ldots, r$. The ordered set of children of a given vertex $v$ of color $i$ may contain:
\begin{itemize}
\item half edges (later called {\it $1$-gons}) linking $v$ to no other vertex. 
\item full edges (later called {\it $2$-gons}) linking $v$ to a vertex of color $i+1$ (modulo $r$). This later vertex is the root of a descending subtree. 
\item $j$-gons linking $v$ to $j-1$ vertices $v_1,v_2,\ldots,v_{j-1}$ of respective colors $i+1,i+2,\ldots,i+j-1$ (modulo $r$). Each $v_k$ is the root of a descending subtree. When the size $j$ of the $j$-gons is not determined, we simply call them {\it polygons}. 
\end{itemize} 

Now we are ready to give the full definition of the considered structure:

\begin{defn}[Cactus trees]\label{defct} For any sequence ${\bf a} = (a_t)$ of $2^r-1$ non-negative integers $a_t$ whith the index $t$ being any non empty strictly increasing subsequence of integers of $1\ldots r$, we define the set $T({\bf a})$ of {\em \bf cactus trees} with vertices of $r$ distinct colors as follows:
\begin{itemize}
\item[(i)] the root vertex of the cactus tree is of color $1$,
\item[(ii)] the ordered set of children of a given vertex $v$ of color $i$ is composed of $j$-gons ($1\leq j \leq r$) linking $v$ to $j-1$ vertices (and subsequent subtrees) of respective colors $i+1,i+2,\ldots,i+j-1$ (modulo $r$), 
\item[(iii)] symbolic labels $\beta_1,\beta_2,\ldots,\beta_n$ (where $n=\sum_t a_t$) are assigned to the polygons such that the set of polygons indexed with the same given label contains exactly one vertex of each color and that all the polygons in the tree are labelled,  
\item[(iv)] for $t=(i_1,i_2,\ldots,i_l)$ ($1\leq i_1< i_2<\ldots<i_l \leq r$), $a_t$ is the number of those sets composed of a $(i_2-i_1)$-gon child of a vertex of color $i_1$, a $(i_3-i_2)$-gon child of a vertex of color $i_2$, $\ldots$, a $(i_l-i_{l-1})$-gon child of a vertex of color $i_{l-1}$ and a $(r-i_l+i_1)$-gon child of a vertex of color $i_l$ with the same symbolic label.\\
(One can easily check that exactly one vertex of each color is contained in such sets.) 
\end{itemize}
\end{defn}

\begin{figure}[htbp]
  \begin{center}
\includegraphics[width=0.60\textwidth]{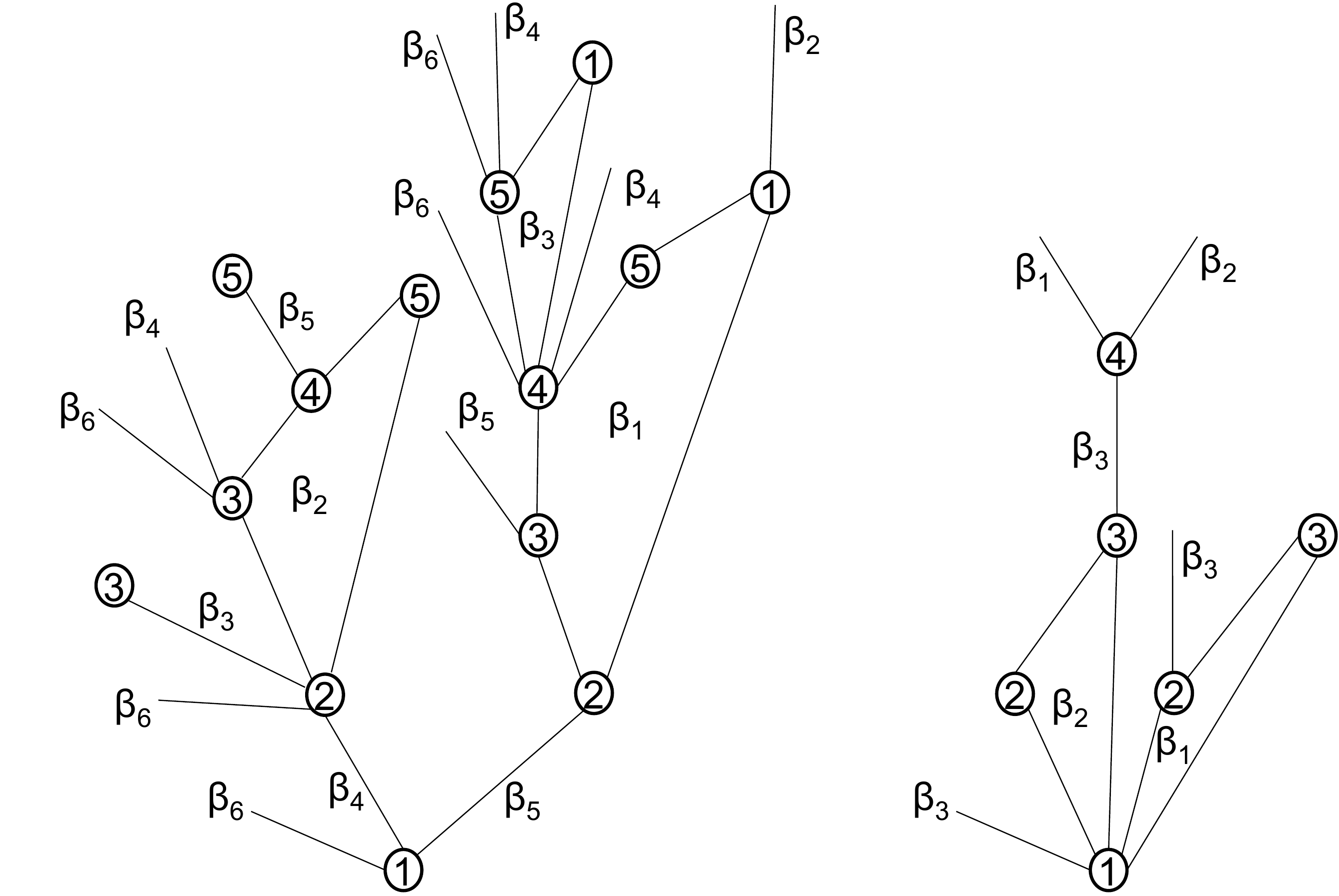}
    \caption{Two examples of cactus trees.}
      \label{CacTreeEx} 
  \end{center}
\end{figure}

\begin{exm}
The cactus tree depicted on the left hand side of Figure \ref{CacTreeEx} has $10$ $1$-gons, $4$ $2$-gons, $1$ $3$-gon, $1$ $4$-gon and $1$ $5$-gon. The corresponding non zero parameters $(a_t)$ are $a_2 = 1$, $a_{1,2} = 1$, $a_{2,4} = 1$, $a_{1,3,4,5} = 1$, $a_{1,3,4} = 1$, $a_{1,2,3,4,5} = 1$. The cactus tree on the right hand side has $4$ $1$-gons, $1$ $2$-gon and $2$ $3$-gons. The $(a_t)$ non equal to zero are $a_{1,4} = 2$, $a_{1,2,3} = 1$.
\end{exm}

\begin{rem}
Point (iii) in Definition \ref{defct} restricts the number of possible cactus trees. In this paper, we consider only the cactus trees for which such a labeling is possible.
\end{rem}

\begin{rem}
If we note $p_i$ the number of vertices of color $i$ ($1\leq i \leq r$)  in a given cactus tree of $T({\bf a})$, we have $\sum_{t;l\notin t}{a_t} = p_l\;$   for $2\leq l \leq r$ and $\sum_{t;1\notin t}{a_t} = p_1-1$.
\end{rem}

\begin{prop} \label{cardcactree}Let ${\bf a} = (a_t)$, $n$, $(p_i)_{1\leq i \leq r}$ be such that $n=\sum_t a_t$ and $\sum_{t;i\notin t}{a_t} = p_i-\delta_{i,1}$. The number $\mid T({\bf a}) \mid$ of cactus trees is given by:
\begin{align}
\label{eqct}
&\mid T({\bf a}) \mid = \frac{{(n-1)!}^{r-1}}{\prod_{1\leq i \leq r}p_i!}\Delta_r({\bf a})\binom{n}{\bf a}
\end{align}
\end{prop}

The proof of Proposition \ref{cardcactree} can be obtained by using the Lagrange theorem in order to compute the number of cactus trees without the labeling and the $1$-gons. Then, counting the number of ways to add the $1$-gons and the symbolic labeling leads to the desired result. Combining Equations (\ref{eqmain}), (\ref{eqlink}) and (\ref{eqct}), we notice that Theorem \ref{thm : main} is equivalent to the following statement:

\begin{thm}
\label{thm : equiv}
The set of partitioned cacti $C^n(p_1,p_2,\ldots,p_r)$ is in bijection with the union of sets of cactus trees $T({\bf a})$ with ${\bf a}$ verifying the properties $n=\sum_t a_t$ and $\sum_{t;i\notin t}{a_t} = p_i-\delta_{i,1}$.
\end{thm}

According to the symmetry property proved in \cite{BM}, Theorem \ref{thm : corro} is implied by Theorem \ref{thm : main}. As a result, Theorem \ref{thm : corro} is also a consequence of Theorem \ref{thm : equiv}. 

\section{Bijection between partitioned cacti and cactus trees}
\label{secbij}
We start with a partitioned cactus $\kappa$ of $C^n(p_1,p_2,\ldots,p_r)$. The edges of the black $r$-gons in $\kappa$ are labeled according to Remark \ref{remcac}. We note $\tilde{\pi}_1,\ldots,\tilde{\pi}_r$ the partitions of the vertices in $\kappa$. Moreover, let $(\alpha_1,\ldots,\alpha_r,\pi_1,\ldots,\pi_r)$ be the $2r$-tuple corresponding to $\kappa$ within the bijection described in Remark \ref{rempartcac}. We proceed with the following construction.\\
First, we define a set containing $p_i$ {\em tree-vertices} (not to be confused with the vertices of $\kappa$) of color $i$ ($1\leq i \leq r$) with an ordered set of children composed of labeled half edges. Each tree-vertex of color $i$ is associated to a block of $\tilde{\pi}_i$ (or equivalently $\pi_i$). The children are indexed by the labels of the edges linking a vertex of color $i$ belonging to the considered block of $\tilde{\pi}_i$ in $\kappa$ and a vertex of color $i-1$ (modulo $r$) {\em sorted in ascending order}. We assume that for all tree-vertices, the resulting labels of the half edges are increasing when we traverse them from {\it left to right}.\\
Then, we look at maximum length sequences of labels $m_i$, $m_{i+1}$, $\ldots$, $m_{i+l}$ (indices are taken modulo $r$) such that : 
\begin{itemize}
\item[(i)] $m_t$ is the greatest label (and therefore the rightmost) of a half edge child of a tree-vertex of color $t$ and 
\item[(ii)] $m_t$'s for $i \leq t \leq i+l$ are also the respective labels of the edges linking the vertices of color $t$ and $t-1$ {\emph{in the same}} $r$-gon of $\kappa$. If such a sequence contains a label $m_1$, maximum index around the tree-vertex of color $1$ that is also parent of a half edge labeled by $1$, we split the sequence into two subsequences $m_i$, $\ldots$, $m_r$ and $m_2$, $\ldots$, $m_{i+l}$. If the initial sequence was the singleton $m_1$, we simply remove it.
\end{itemize}
\begin{lem}
The maximum length of the sequences defined above is $r-1$.
\end{lem}
Next we build a cactus tree by connecting the tree-vertices using $j$-gons ($j\geq 2$). The tree-vertex of color $1$ connected to the half-edge labeled $1$ is the root of the cactus tree. For any sequence $m_i$, $m_{i+1}$, $\ldots$, $m_{i+l}$ corresponding to the same $r$-gon in $\kappa$, let $t$ be the label of the edge linking the vertex of colors $i-1$ and $i-2$ in this $r$-gon. By definition, $t$ is not the maximum label around a tree-vertex of color $i-1$. We connect the tree-vertices with maximum child labels $m_i$, $m_{i+1}$, $\ldots$, $m_{i+l}$ thanks to a $l+2$-gon. We substitute this $l+2$-gon to the half edge labeled $t$ in the children set of the corresponding tree-vertex of color $i-1$. We assign the labels $m_i$, $m_{i+1}$, $\ldots$, $m_{i+l}$, $t$ to it so that the edge linking the tree-vertices of colors $i$ and $i-1$ is labeled $m_i$, the edge linking the tree-vertices of colors $i+1$ and $i$ is labeled $m_{i+1}$, and so on, the edge linking the tree-vertices of colors $i+l$ and $i+l-1$ is labeled $m_{i+l}$ and the edge linking the tree-vertices of colors $i+l$ and $i-1$ is labeled $t$. In what follows, we use the term  $u$-{\em color numeric label} of the considered $l+2$-gon for $m_u$ ($i\leq u \leq i+l$) and $i-1$-{\em color numeric label} for $t$. \\
Finally, we allocate a symbolic label from $\beta_1,\ldots,\beta_n$ to each polygon and each of the remaining half edges ($1$-gons) such that all the $j$-gons ($1\leq j \leq r$) with numeric labels corresponding to the same $r$-gon of $\kappa$ have the same symbolic label. At this stage, we remove all the numeric labels.

\begin{figure}[htbp]
  \begin{center}
	\includegraphics[width=0.60\textwidth]{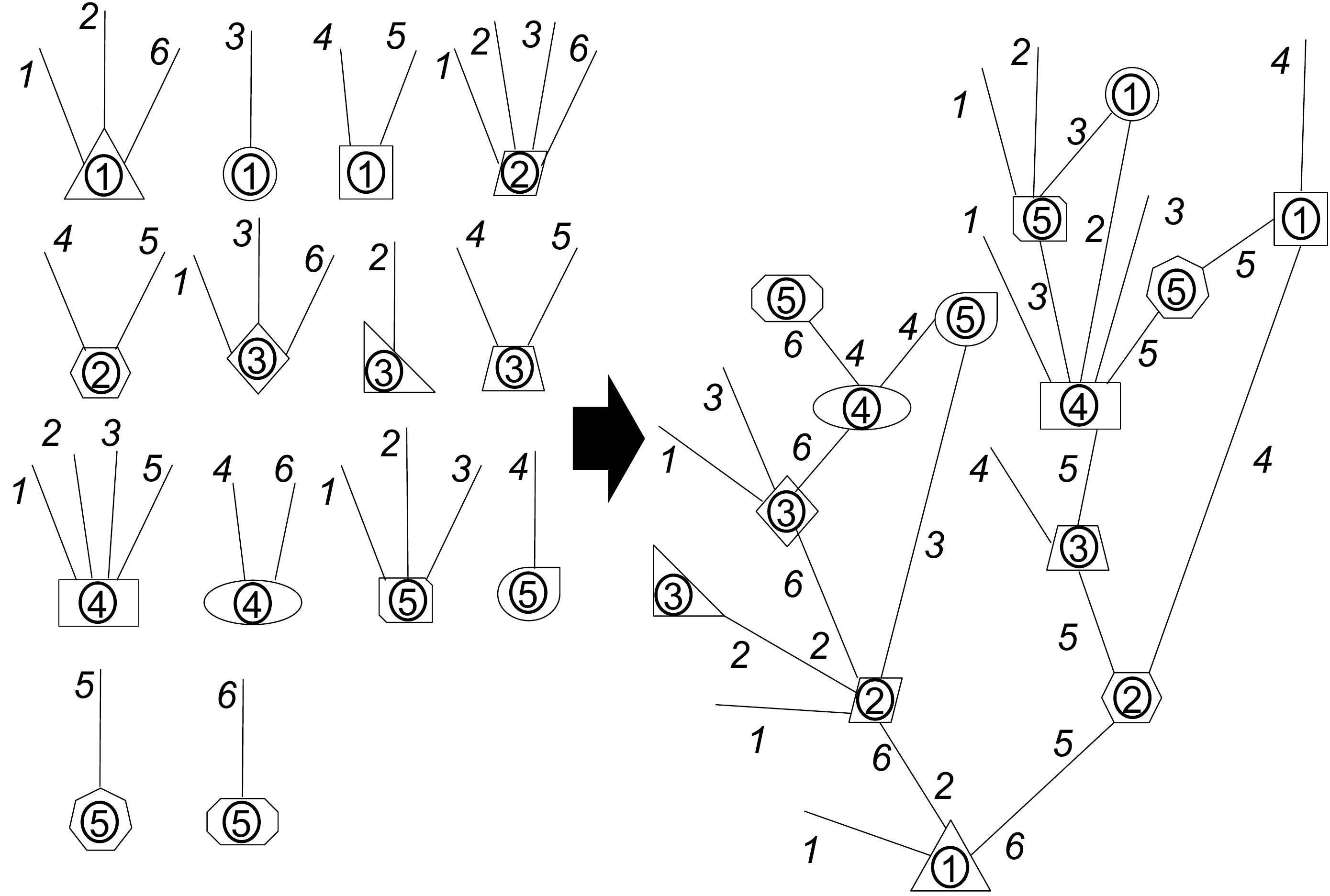}
    \caption{Application of the bijective construction to the partitioned cactus of Figure \ref{cactipart}}
    \label{figbij}
  \end{center}
\end{figure}

\begin{exm} We apply the construction above to the partitioned cactus depicted on Figure \ref{cactipart}. The definition of the set of tree-vertices and their connections by polygons is shown on Figure \ref{figbij}. The final step of symbolic labeling leads to the cactus on the left hand side of Figure \ref{CacTreeEx}. 
\end{exm}

\begin{lem}
The construction above defines a cactus tree $\tau$ in $T({\bf a})$ for some ${\bf a}$ verifying the properties $n=\sum_t a_t$ and $\sum_{t;i\notin t}{a_t} = p_i-\delta_{i,1}$. 
\end{lem}

\section{Proof of the bijection}
\label{secproof}
\subsection{Injectivity}
\label{secinj}
Assume $\tau$ is a cactus tree of $T({\bf a})$ obtained by the construction of Section \ref{secbij}. We show that at most one $2r$-tuple $(\alpha_1,\ldots,\alpha_r,\pi_1,\ldots,\pi_r)$ (or equivalently at most one partitioned cactus) is mapped to $\tau$. To this purpose we show by induction that the numeric labels removed at the end of the procedure can be uniquely recovered.\\
First step is to notice that $1$ is necessarily (one of) the numeric label(s) of the leftmost child of the root of $\tau$ (of color $1$).\\
Now assume that the $u$-color (resp. $1$-color) numeric labels $1,\ldots, i-1$ (resp. $1,\ldots, i$) have been recovered for $u=2,3,\ldots,r$ and $i<n$.  
\begin{itemize}
\item Let $\beta$ be the symbolic label of the polygon in $\tau$ with recovered numeric $1$-color label $i$ (this polygon is by assumption connected to a vertex of color $1$). Then, $\beta$ is also the symbolic label of exactly one polygon (possibly the one with recovered numeric $1$-color label $i$) connected to a vertex $v$ of color $r$. But as noticed in Lemma \ref{lemcac}, if the $1$-color numeric label is $i$, the $r$-color numeric label corresponding to the same $r$-gon in the initial partitioned cactus is $\alpha_r^{-1}(i)$. As a result, $\alpha_r^{-1}(i)$ is the numeric $r$-color label of the polygon connected to $v$. The blocks of $\pi_r$ are stable by $\alpha_r$ so $\alpha_r(i)$ belongs to the same block of $\pi_r$. According to the order of the $r$-color labels around $v$, $i=\alpha_r^{-1}\alpha_r(i)$ is necessarily the $r$-color numeric label of the leftmost polygon connected to $v$ with non recovered $r$-color label. 
\item Assume that we have recovered the polygon with $u$-color label $i$ ($2 < u$). The index of the corresponding black $r$-gon in the partitioned cactus is necessarily $\alpha_u\alpha_{u+1}\ldots\alpha_r(i)$. Integers $\alpha_u\alpha_{u+1}\ldots\alpha_r(i)$ and $\alpha_{u-1}\alpha_u\alpha_{u+1}\ldots\alpha_r(i)$ belong to the same block of $\pi_{u-1}$. As the $(u-1)$-color labels have been sorted according to $\alpha_{r}^{-1}\ldots\alpha_{u-1}^{-1}$ around the vertices of color $(u-1)$, $i$ is necessarily the $(u-1)$-color numeric label of the leftmost polygon (for which such a label has not been recovered yet) around the vertex of color $u-1$ incident to the polygon with the same symbolic label as the one incident to the polygon of $u$-color label $i$. 
\item Finally assume that we have recovered the polygon with $2$-color label $i$. The index of the corresponding black $r$-gon in the partitioned cactus is necessarily $\alpha_2\alpha_{3}\ldots\alpha_r(i)$. We use the symbolic label to identify the vertex $v$ of color $1$ incident to the polygon of $1$-color label $\alpha_2\alpha_{3}\ldots\alpha_r(i)$. With a similar argument as above $i+1= \alpha_1\alpha_2\alpha_{3}\ldots\alpha_r(i)$ is necessarily the $1$-color label of the leftmost polygon (with no recovered $1$-color numeric label) incident to $v$.
\end{itemize}
The knowledge of $\tau$ uniquely determines the numeric labels of the polygons. But it is easy to see that combining symbolic and numeric labels uniquely determines $\alpha^{-1}_r,\ldots,\alpha_{r}^{-1}\ldots\alpha_{u}^{-1},\ldots,\alpha_{r}^{-1}\ldots\alpha_{2}^{-1}$ and the $\alpha_i$ themselves. Then the knowledge of the numeric labels around the same vertices of $\tau$ uniquely determines the partitions $\pi_1,\ldots,\pi_r$. As a result, the partitioned cactus is uniquely determined.

\begin{figure}[htbp]
  \begin{center}
   \includegraphics[width=0.60\textwidth]{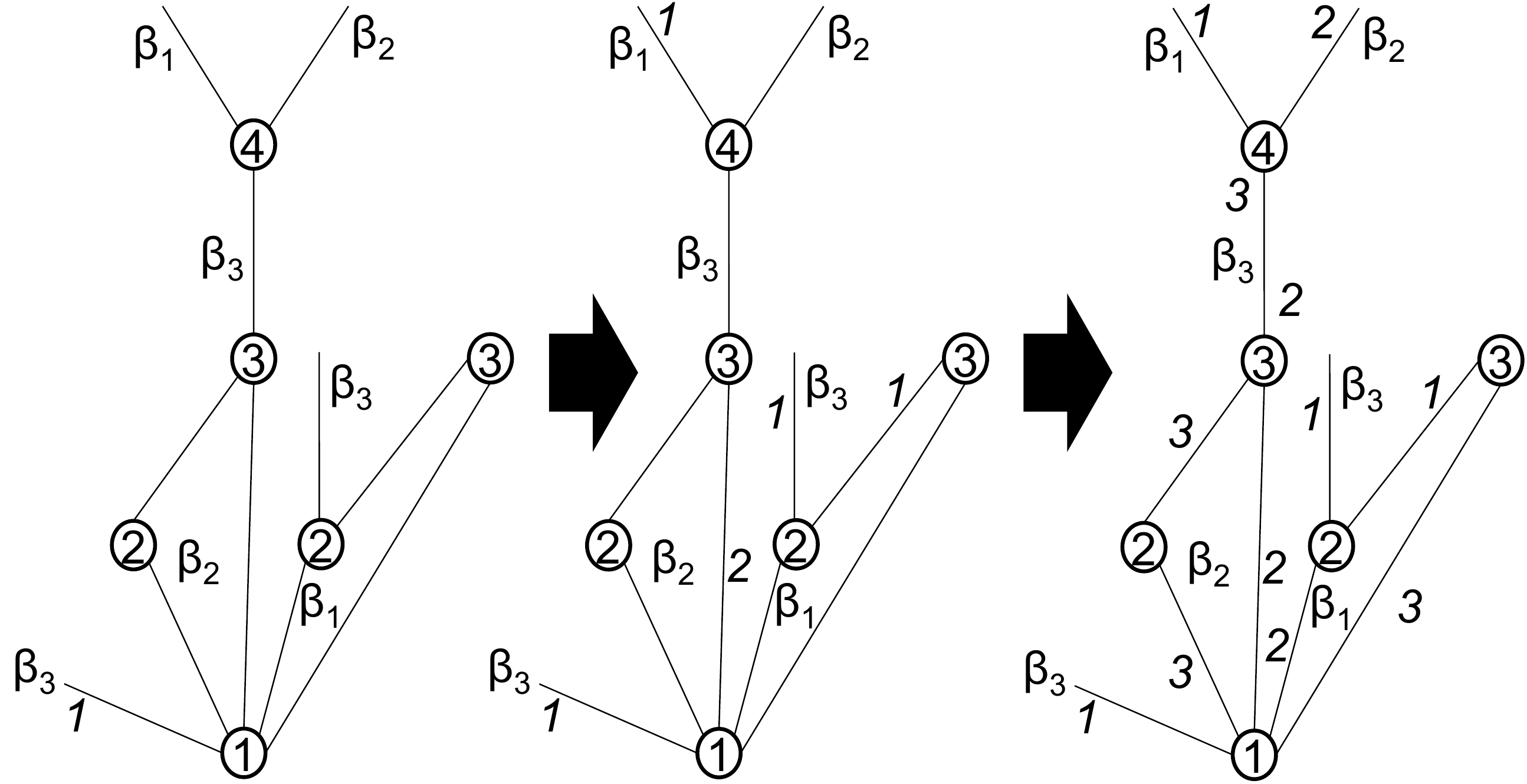}
    \caption{Application of the inverse procedure to the cactus tree on the right hand side of Figure \ref{CacTreeEx}.}
      \label{figrecons} 
  \end{center}
\end{figure}

\begin{exm}
We apply this inverse procedure to the cactus tree on the right hand side of Figure \ref{CacTreeEx}. Figure \ref{figrecons} shows how the numeric labels are iteratively recovered. The resulting numerically labeled cactus tree corresponds to the cactus on the right hand side of Figure \ref{cacti} with partitions $\pi_1 = \{1,2,3\}$, $\pi_2 = \{1,3\}\{2\}$, $\pi_3 = \{1,2\}\{3\}$ and $\pi_4 = \{1,2,3\}$.
\end{exm}

\subsection{Surjectivity}
We show that the reconstruction procedure defined in Section \ref{secinj} always ends with a valid output. The procedure would be interrupted before the full recovery of the numeric labels if and only if no leftmost polygon around a vertex $v$ of color $i$ with non recovered $i$-color numeric label is available. Two cases are to be considered.
\begin{itemize}
\item[(i)] If $v$ is not the root, this situation is clearly impossible. If the number of polygons incident to $v$ is $c$, we traverse exactly $c$ times the vertex $v$ to allocate $i$-color labels (the symbolic labels link the polygons around $v$ to exactly $c$ polygons incident to vertices of color $i+1$).
\item[(ii)] If $v$ is the root, then an additional difficulty occurs as $1$-color label $1$ is recovered out of the main procedure. However for any vertex $u$ of color $j+1$ in the cactus tree, the rightmost polygon connecting it to a vertex $w$ of color $j$ is obviously the last one to be recovered. The symbolic label of this polygon naturally links $u$ and $w$. As a result, all the $j$-color labels of $w$ are recovered necessarily {\it after} all the $j+1$ - labels are recovered around $u$. This extends to all the children of $w$ and the vertex of color $j+1$ in the rightmost polygon of $w$ if any. Since $v$ is the root all the $1$-color labels of $v$ are recovered necessarily after all the labels of all the other vertices in the cactus tree are recovered. Finally the procedure ends in the proper way.
\end{itemize}

\section{Equivalence of the main theorem and Jackson's formula}
\label{seceq}
A natural question is the equivalence between the formula of Theorem \ref{thm : main} and Jackson's formula of \cite{DMJ} addressing the factorizations of $\gamma_n$. 
The result of Jackson for the factorizations of the long cycle can be stated as follows:
\begin{align}
\nonumber \frac{1} {(n!)^{r-1}}&\sum_{p_1, p_2,\ldots, p_r} {k^n_{p_1, p_2,\ldots, p_r}} \prod_{1\leq i \leq r}x_i^{p_i}\\ 
&=  \phi\left (\prod_{1\leq i \leq r}x_i\left(\prod_{1\leq i \leq r}(1+x_i)-\prod_{1\leq i \leq r}(x_i)\right)^{n-1}\right)
\end{align}
\noindent where $\phi$ is the mapping defined by $\phi(\prod_{i}x_i^{p_i})=\prod_{i}\binom{x_i}{p_i}$ and extended linearly. One can show that this formula is equivalent to:
\begin{align}
\frac{1} {(n!)^{r-1}}&\sum_{p_1, p_2,\ldots, p_r} {k^n_{p_1, p_2,\ldots, p_r}} \prod_{1\leq i \leq r}x_i^{p_i} =  \sum_{p_1, p_2,\ldots, p_r}\sum_{{\bf a}}\binom{n-1}{\bf a}\prod_{1\leq i \leq r}\binom{x_i}{p_i}
\end{align}
\noindent where the sequences ${\bf a} = (a_t)$ of $2^r-1$ non-negative integers $(a_t)$ satisfy
$
\sum_t{a_t} = n-1,\;\;\sum_{t;l\in t}{a_t} = p_l-1 \;\mbox{ for }\; 1\leq l \leq r
$.\\ 
For $r=2$ the equivalence is obvious as in Equation (\ref{eqmain}) the only sequence ${\bf a}$ fitting the conditions is $a_{1} = p_2$, $a_{2} = p_1-1$, $a_{1,2} = n+1-p_1-p_2$ and $\Delta_2({\bf a}) = p_2$ in this case. The summand for indices $p_1$ and $p_2$ in our formula reads
\begin{align}
(n-1)!\,p_2\binom{n}{p_1-1,p_2}=n!\binom{n-1}{p_1-1,p_2-1}
\end{align}
This shows that the two results are identical in this case.\\ 
For $r=3$, the determinant is $\Delta_3({\bf a}) = p_2p_3-a_3(p_3-a_1)$.
For given $p_1$, $p_2$ and $p_3$ the equivalence between the two formulas can be shown in two steps. First we have:
\footnotesize
\begin{align}
\nonumber&(p_2p_3-a_3(p_3-a_1))\binom{n}{a_1, a_2, a_3, p_1-1-a_2-a_3,p_2-a_3-a_1,p_3-a_1-a_2}\\
\nonumber &=((p_2-a_3-a_1)(p_3-a_2-a_1)+a_1(p_2+p_3-a_1)+a_2(p_2-a_3-a_1))\\
\nonumber&\;\;\;\;\;\;\;\;\;\times \binom{n}{a_1, a_2, a_3, p_1-1-a_2-a_3,p_2-a_3-a_1,p_3-a_1-a_2}\\
\nonumber&=n(n+2+a_1+a_2+a_3 -p_1-p_2-p_3)\\
\nonumber&\;\;\;\;\;\;\;\;\;\times\binom{n-1}{a_1, a_2, a_3, p_1-1-a_2-a_3,p_2-1-a_3-a_1,p_3-1-a_1-a_2}\\
\nonumber&\;\;\;\;\;\;\;\;\;+n(p_1-a_3-a_2)\binom{n-1}{a_1, a_2-1, a_3, p_1-a_2-a_3,p_2-1-a_3-a_1,p_3-a_1-a_2}\\
\nonumber&\;\;\;\;\;\;\;\;\;+n(p_2+p_3-a_1)\binom{n-1}{a_1-1, a_2, a_3, p_1-1-a_2-a_3,p_2-a_3-a_1,p_3-a_1-a_2}
\end{align}
\normalsize
Then summing over $a_1$, $a_2$, $a_3$ with the proper shifts of variable to get the same multinomial coefficient brings us to
\footnotesize
\begin{align}
\nonumber&(n-1)!^2\sum_{a_1,a_2,a_3}(p_2p_3-a_3(p_3-a_1))\binom{n}{a_1, a_2, a_3, p_1-1-a_2-a_3,p_2-a_3-a_1,p_3-a_1-a_2}\\
\nonumber &=n!^2\sum_{a_1,a_2,a_3}\binom{n-1}{a_1, a_2, a_3, p_1-1-a_2-a_3,p_2-1-a_3-a_1,p_3-1-a_1-a_2}
\end{align}
\normalsize
that proves the equivalence of the two formulas also in the case $r=3$.


\bibliographystyle{alpha}

%

%
%
%

\end{document}